\def\hldest#1#2#3{}


\font\smallheader=cmssbx10 

\font\eightpt=cmr8
\font\ninept=cmr9


\font\ninerm=cmr9     \font\eightrm=cmr8   \font\sixrm=cmr6      
\font\ninei=cmmi9     \font\eighti=cmmi8   \font\sixi=cmmi6      
\font\ninesy=cmsy9    \font\eightsy=cmsy8  \font\sixsy=cmsy6     
\font\ninebf=cmbx9    \font\eightbf=cmbx8  \font\sixbf=cmbx6     
\font\ninett=cmtt9    \font\eighttt=cmtt8                        
\font\nineit=cmti9    \font\eightit=cmti8     
\font\ninesl=cmsl9    \font\eightsl=cmsl8                        

\font\tensc=cmcsc10   \font\ninesc=cmcsc9  \font\eightsc=cmcsc8  

\font\eightssq=cmssq8  \font\eightssqi=cmssqi8  


\font\tenssbx=cmssbx10 

  \font\twelvebf=cmbx12
  
\def\sc{\tensc}  \def\mc{\ninerm}


\newskip\ttglue
\def\tenpoint{\def\rm{\fam0\tenrm}%
  \textfont0=\tenrm \scriptfont0=\sevenrm \scriptscriptfont0=\fiverm
  \textfont1=\teni  \scriptfont1=\seveni  \scriptscriptfont1=\fivei
  \textfont2=\tensy \scriptfont2=\sevensy \scriptscriptfont2=\fivesy
  \textfont3=\tenex \scriptfont3=\tenex   \scriptscriptfont3=\tenex
  \textfont\itfam=\tenit  \def\it{\fam\itfam\tenit}%
  \textfont\slfam=\tensl  \def\sl{\fam\slfam\tensl}%
  \textfont\ttfam=\tentt  \def\tt{\fam\ttfam\tentt}%
  \textfont\bffam=\tenbf  \scriptfont\bffam=\sevenbf
   \scriptscriptfont\bffam=\fivebf \def\bf{\fam\bffam\tenbf}%
  \tt \ttglue=.5em plus.25em minus.15em
  \normalbaselineskip=12pt
  \setbox\strutbox=\hbox{\vrule height8.5pt depth3.5pt width0pt}%
  \let\sc=\tensc \let\mc=\ninerm  
  \def\cyr{\tencyr\cyracc}\def\cyri{\tencyri\cyracc}
  \let\big=\tenbig  \normalbaselines\rm}

\def\ninepoint{\def\rm{\fam0\ninerm}%
\textfont0=\ninerm  \scriptfont0=\sixrm  \scriptscriptfont0=\fiverm
\textfont1=\ninei   \scriptfont1=\sixi   \scriptscriptfont1=\fivei
\textfont2=\ninesy  \scriptfont2=\sixsy  \scriptscriptfont2=\fivesy
\textfont3=\tenex   \scriptfont3=\tenex  \scriptscriptfont3=\tenex
\textfont\itfam=\nineit  \def\it{\fam\itfam\nineit}%
\textfont\slfam=\ninesl  \def\sl{\fam\slfam\ninesl}%
\textfont\ttfam=\ninett  \def\tt{\fam\ttfam\ninett}%
\textfont\bffam=\ninebf  \scriptfont\bffam=\sixbf
\scriptscriptfont\bffam=\fivebf\def\bf{\fam\bffam\ninebf}%
\tt\ttglue=.5em plus.25em minus.15em
\normalbaselineskip=11pt
\setbox\strutbox=\hbox{\vrule height8pt depth3pt width0pt}%
\let\sc=\ninesc\let\mc=\eightrm
\def\cyr{\ninecyr\cyracc}\def\cyri{\ninecyri\cyracc}
\let\big=\ninebig\normalbaselines\rm}

\def\eightpoint{\def\rm{\fam0\eightrm}%
  \textfont0=\eightrm \scriptfont0=\sixrm \scriptscriptfont0=\fiverm
  \textfont1=\eighti  \scriptfont1=\sixi  \scriptscriptfont1=\fivei
  \textfont2=\eightsy \scriptfont2=\sixsy \scriptscriptfont2=\fivesy
  \textfont3=\tenex   \scriptfont3=\tenex \scriptscriptfont3=\tenex
  \textfont\itfam=\eightit  \def\it{\fam\itfam\eightit}%
  \textfont\slfam=\eightsl  \def\sl{\fam\slfam\eightsl}%
  \textfont\ttfam=\eighttt  \def\tt{\fam\ttfam\eighttt}%
  \textfont\bffam=\eightbf  \scriptfont\bffam=\sixbf
  \normalbaselineskip=9pt
  \let\sc=\eightsc \let\mc=\sevenrm  
  \def\cyr{\eightcyr\cyracc}\def\cyri{\eightcyri\cyracc}
  \let\big=\eightbig  \normalbaselines\rm}%
\def\nospace{\nulldelimiterspace0pt\mathsurround0pt}%
\def\tenbig#1{{\hbox{$\left#1\vbox to8.5pt{}\right.\nospace$}}}%
\def\ninebig#1{{\hbox{$\textfont0=\tenrm\textfont2=\tensy
  \left#1\vbox to7.25pt{}\right.\nospace$}}}%
\def\eightbig#1{{\hbox{$\textfont0=\ninerm\textfont2=\ninesy
  \left#1\vbox to6.5pt{}\right.\nospace$}}}%

\def\nonextendedbold{
  \font\fiveb=cmb10 at 5pt
  \font\sixb=cmb10 at 6pt
  \font\sevenb=cmb10 at 7pt
  \font\eightb=cmb10 at 8pt
  \font\nineb=cmb10 at 9pt
  \font\tenb=cmb10
  \font\twelveb=cmb10 at 12pt
  \let\fivebf=\fiveb
  \let\sixbf=\sixb
  \let\sevenbf=\sevenb
  \let\eightbf=\eightb
  \let\ninebf=\nineb
  \let\tenbf=\tenb
  \let\twelvebf=\twelveb
}

\def\leftrighttop#1#2{
  \headline{\ifnum\pageno=1\hfil\else{\ninept #1 \hfil #2}\fi}
}

\def\firstnopagenum{
  \footline{\ifnum\pageno=1 \hfil \else \hfil{\rm \number\pageno}\hfil\fi}
}

\def\maketitle#1#2#3#4{
  \centerline {\titlefont #1}
  \medskip
  \centerline {\eightpt #2}
  \medskip
  \centerline {\tensc #3}
  \medskip
  \centerline {\tensc #4}
  \bigskip
}


\outer\def\floattext#1 #2. #3\par{
  $$
  \vbox{
    \hsize #1 true in
    \noindent{\bf #2.}\enskip #3
  }
  $$
}


\def\lsection#1\par{
  \bigskip\vskip\parskip
  \leftline{\sectionfont#1}\nobreak\medskip\noindent
}

\def\csection#1\par{
  \bigskip\vskip\parskip
  \centerline{\sectionfont#1}\nobreak\medskip\noindent
}

\def\rsection#1\par{
  \bigskip\vskip\parskip
  \rightline{\sectionfont#1}\nobreak\medskip\noindent
}
\def\section{\lsection}

\def\boldlabel#1. {\noindent{\bf #1.\enspace}}
\def\subsection#1. {\medskip\noindent{\bf #1.\enspace}}


\def\begincode{\par\begingroup\nobreak\medskip \obeylines\obeyspaces\frenchspacing\tt}
\def\endcode{\endgroup\medbreak\noindent}

\font\tenfrak=eufm10
\font\sevenfrak=eufm7
\font\fivefrak=eufm5
\newfam\frakfam
\textfont\frakfam=\tenfrak
\scriptfont\frakfam=\sevenfrak
\scriptscriptfont\frakfam=\fivefrak

\def\janksc#1#2 {#1{\eightpt#2}}
\def\jankscsp#1#2 {#1{\eightpt#2}\ }
\def\scproclaim#1.#2\par{\noindent\jankscsp #1.\enspace{\it#2\par}}

\def\bref#1{[#1]}
\def\ref#1{[#1]}

\def\quote{
  \begingroup
    \baselineskip 10pt
    \parfillskip 0pt
    \interlinepenalty 10000
    \leftskip 0pt plus 40pc minus \parindent
    \let\rm=\quoterm\let\sl=\quotesl\everypar{\sl}
    \obeylines
}
\def\author#1(#2){\nobreak\smallskip\rm--- \rm#1\unskip\enspace(#2)\par\endgroup}

\def\titlefont{\twelvebf}
\def\sectionfont{\tenssbx}
\def\quoterm{\eightssq}
\def\quotesl{\eightssqi}



\def\bookheader#1#2{
  \nopagenumbers
  \def\leftheadline{{\rm\folio}\hfil{\eightpoint#1}\hfil}
  \def\rightheadline{\hfil{\eightpoint#2}\hfil{\rm\folio}}
  \headline{\ifodd\pageno{\ifnum\pageno<2\hfil\else\rightheadline\fi}\else\leftheadline\fi}
}

\tenpoint



\def\xskip{\hskip 7pt plus 3pt minus 4pt}

\def\proof{\medbreak\noindent{\it Proof.}\xskip\ignorespaces}

\def\slug{\quad\hbox{\kern1.5pt\vrule width2.5pt height6pt depth1.5pt\kern1.5pt}\medskip}
\def\noskipslug{\quad\hbox{\kern1.5pt\vrule width2.5pt height6pt depth1.5pt\kern1.5pt}}

\newdimen\algindent
\newif\ifitempar \itempartrue 
\def\algindentset#1{\setbox0\hbox{{\bf #1.\kern.25em}}\algindent=\wd0\relax}
\def\algbegin #1 #2{\algindentset{#21}\alg #1 #2} 
\def\aalgbegin #1 #2{\algindentset{#211}\alg #1 #2} 
\def\alg#1(#2). {\medbreak 
  \noindent{\bf#1}({\it#2\/}).\xskip\ignorespaces}
\def\algstep#1.{\ifitempar\smallskip\noindent\else\itempartrue
  \hskip-\parindent\fi
  \hbox to\algindent{\bf\hfil #1.\kern.25em}%
  \hangindent=\algindent\hangafter=1\ignorespaces}


\def\NN{{\bf N}}
\def\ZZ{{\bf Z}}


\def\op#1{\mathop{\hbox{#1}}\nolimits}


\def\one{\op{\bf 1}}



\def\eps{\epsilon}

\newcount\thmcount  
\thmcount=1
\newcount\sectcount  
\sectcount=1
\newcount\figcount  
\figcount=1
\newcount\eqcount  
\eqcount=1

\def\oldno#1{\eqno({\oldstyle#1})}

\def\adveq{\oldno{\the\eqcount}\global\advance\eqcount by 1}  
\def\advthm{\the\thmcount\global\advance \thmcount by 1}

\def\advsect{\section\the\sectcount\global\advance\sectcount by 1. }

\def\caption#1{\centerline{\ninepoint{\bf Fig.~\the\figcount\global\advance\figcount by 1.\enspace}#1}}

\outer\def\parenproclaim #1 (#2).#3\par{\medbreak
  \noindent{\bf #1}\enspace\rm({\it #2\/}).\nobreak\ignorespaces{\sl #3\par}
  \ifdim\lastskip<\medskipamount \removelastskip\penalty55\medskip\fi}


\newdimen\axiomindent
\def\axiomindentset#1{\setbox0\hbox{{\bf #1.\kern.25em}}\axiomindent=\wd0\relax}
\def\axiom#1. [#2.]{\ifitempar\par\noindent\else\itempartrue
  \hskip-\parindent\fi%
  \hbox to\axiomindent{\bf\hfil #1.\kern.25em}%
  \hangindent=\axiomindent\hangafter=1[{\it #2.}]}

\input epsf
\input eplain

\def\top{\widehat 1}
\def\bottom{\widehat 0}
\def\rank{\op{\rm rank}}

\def\ref#1{\special{ps:[/pdfm { /big_fat_array exch def big_fat_array 1 get 0
0 put big_fat_array 1 get 1 0 put big_fat_array 1 get 2 0 put big_fat_array pdfmnew } def}%
[\hlstart{name}{}{bib#1}#1\hlend]}

\def\tilde{\widetilde}

\def\F{{\cal F}}
\def\tab{\ \ \ \ }

\magnification=\magstephalf
\hoffset=40pt \voffset=28pt
\hsize=29pc  \vsize=45pc  \maxdepth=2.2pt  \parindent=19pt
\nopagenumbers
\def\leftheadline{{\rm\folio}\hfil{\eightpoint GOH AND SAKS}\hfil}
\def\rightheadline{\hfil{\eightpoint ON THE HOMOLOGY OF SEVERAL NUMBER-THEORETIC SET FAMILIES}\hfil{\rm\folio}}
\headline={\ifodd\pageno{\ifnum\pageno<2\hfil\else\rightheadline\fi}\else\leftheadline\fi}

\enablehyperlinks

\ifpdf
  \hlopts{bwidth=0}
  \pdfoutline goto name {intro} {Introduction}%
  \pdfoutline goto name {proof} {The Mobius function and alternating sums}%
  \pdfoutline goto name {prim} {Primitive, pairwise coprime, and product-free sets}%
  \pdfoutline goto name {coprimefree} {Coprime-free sets}%
  \pdfoutline goto name {general} {A generalization of primitive sets}%
  \pdfoutline goto name {acks} {Acknowledgements}%
  \pdfoutline goto name {refs} {References}%
  \pdfoutline goto name {numerical} {Numerical tables}%
  \pdfoutline goto name {verification} {Computer verification of Proposition 8}%
\else
  \special{ps:[/PageMode /UseOutlines /DOCVIEW pdfmark}%
  \special{ps:[/Count -0 /Dest (intro) cvn /Title (Introduction) /OUT pdfmark}%
  \special{ps:[/Count -0 /Dest (proof) cvn /Title (The Mobius function and alternating sums) /OUT pdfmark}%
  \special{ps:[/Count -0 /Dest (prim) cvn
               /Title (Primitive, pairwise coprime, and product-free sets) /OUT pdfmark}%
  \special{ps:[/Count -0 /Dest (coprimefree) cvn
               /Title (Coprime-free sets) /OUT pdfmark}%
  \special{ps:[/Count -0 /Dest (general) cvn /Title (A generalisation of primitive sets) /OUT pdfmark}%
  \special{ps:[/Count -0 /Dest (acks) cvn /Title (Acknowledgements) /OUT pdfmark}%
  \special{ps:[/Count -0 /Dest (refs) cvn /Title (References) /OUT pdfmark}%
  \special{ps:[/Count -0 /Dest (numerical) cvn /Title (Numerical tables) /OUT pdfmark}%
  \special{ps:[/Count -0 /Dest (verification) cvn /Title (Computer verification of Proposition 8) /OUT pdfmark}%
\fi

\centerline{\twelvebf On the homology of several number-theoretic set families}
\bigskip
\bigskip

\centerline{\sc Marcel K. Goh}
\smallskip
\centerline{\sl Edmonton, Alta., Canada}
\smallskip
\centerline{\tt marcel.goh@mail.mcgill.ca}

\bigskip
\medskip
\centerline{\sc Jonah Saks}
\smallskip
\centerline{\sl Department of Mathematics and Statistics, McGill University}
\smallskip
\centerline{\sl Montreal, Que., Canada}
\smallskip
\centerline{\tt jonah.saks@mail.mcgill.ca}

\floattext5 \ninebf Abstract.
\ninepoint
This paper describes the homology of various simplicial complexes associated to set families from combinatorial
number theory, including primitive sets, pairwise coprime sets, product-free sets, and coprime-free sets.
We present a condition on a set family that results in easy computation of the homology groups, and show that the first
three examples, among many others, admit such a structure. We then extend our techniques
to address the complexes associated to coprime-free sets, as well as a generalisation of primitive sets.
\smallskip
\noindent\boldlabel Keywords. Homology, primitive sets, product-free sets, coprime-free sets.
\smallskip
\noindent\boldlabel MSC classes. 06A07, 11B75.

\advsect Introduction
\hldest{xyz}{}{intro}

{\sc Homology theory} has been a fundamental tool in enumerative combinatorics since
G.-C.~Rota's watershed 1964 paper~\bref{14} established links between the topological data of a poset
and its combinatorial properties. Combinatorial problems, reciprocally, give rise to
a wealth of topological spaces to study. In particular, many families of sets that arise in
the field of combinatorial number theory possess interesting topological structures.

In this paper we consider a number of simplicial
complexes associated to set families defined by conditions of a number-theoretic nature. It turns out that
many of these complexes are in some sense the same, topologically speaking, and we introduce a definition---namely,
that of a partition-intersecting family---to capture a condition that engenders this structure.
The remainder of the introduction will be devoted to this task. In the second section we detail well-known
relationships between the homology groups of a poset and alternating sums associated to its elements.
Three that fall into our definition of
partition-intersecting are primitive sets, pairwise coprime sets, and product-free sets.
In Section~3 we describe these homologies without
breaking a sweat, and enumerate a few more examples that can be dealt with in the same way.
On the other hand, the family of coprime-free sets is not partition-intersecting, and in this case we can only
fully characterize the zeroth and first homology groups. This is done in Section~4.
In Section~5, we tackle an analogous generalisation of primitive sets,
whose corresponding homology groups are
derivable nearly free of charge from the machinery used to deal with coprime-free sets.

The number-theoretic properties of the set families underlying this paper have been extensively studied; we have
chosen to defer the exposition of this background to Section~3, when these sets are defined in more detail.

\medskip\boldlabel Definitions and notation.
For integers $n\ge 1$ we denote by $[n]$ the discrete interval $\{1,\ldots,n\}$. The main object of study will
be some family $\F$ of finite subsets of integers, so for brevity we shall let
$\F_n$ denote the finite set $\F\cap 2^{[n]}$.
Consider the partially ordered set (poset) obtained by taking $\F_n$ as a ground set and ordering
its elements by inclusion. If $\emptyset\in \F_n$, then it is the bottom element of the poset,
since $\emptyset \subseteq x$
for all $x\in \F_n$. It will become important later to add the restriction $[n]\notin \F_n$ for all $n\ge 2$,
so that $\F_n$ does not have a top element, but in this case we may adjoin an artificial top
element to obtain a lattice; we shall denote this element by $\top$.
An element $y$ is said to {\it cover} an element $x\ne y$ if $x\le y$ and $x\le z\le y$ implies that
$z = x$ or $z=y$. In a lattice, an element covered by $\top$ is called a {\it coatom}.

\medskip\boldlabel Cross-cuts. A {\it chain} in a poset $X$ is a set $\{x_1,x_2,\ldots,x_k\}\subseteq X$
such that $x_1 < x_2 < \cdots < x_k$. A chain is {\it maximal}\/ if adding any new element to the chain
causes the chain property to be violated.
Let $L$ be a lattice with top element $\top$ and bottom element $\bottom$.
A {\it cross-cut} $C$ is a subset of $L\setminus\{\top,\bottom\}$ such that no two elements of $C$ are
comparable and every chain in $L$ contains some element of $C$.
It is easy to see that the set of coatoms in a
lattice is a cross-cut.

A subset $S\subseteq L$ is said to be {\it spanning} if the join of all its elements is $\top$ and
the meet of all its elements is $\bottom$. The cross-cut complex $\Delta(C)$ of a cross-cut
$C$ is the simplicial complex obtained from taking $C$ as the vertex set and adding a face for any subset of
$C$ that is {\it not} spanning.
(For us, a {\it simplicial complex} is a set of subsets of a finite set that is closed under intersection.)
It is a classical result that the homology groups of $\Delta(C)$ depend only on the lattice
$L$ and not the choice of cross-cut $C$, so it makes sense to speak of the {\it homology} of $L$.
For the rest of the paper, by homology we shall
mean reduced homology (i.e., the zeroth homology of a connected
space is trivial) over the group $\ZZ$. So we write $\tilde H_k(\Delta)$ for $\tilde H_k(\Delta,\ZZ)$.

\medskip\boldlabel Partition-intersecting families.
We now define a criterion on a poset that leads to easily-described homology groups.
Let $\F$ be a family of subsets of $\NN$. For $n\ge 2$ in the definition below, let
$\F_n = \F\cap 2^{[n]}$ for brevity and let $M_n$ be the set of maximal elements in $\F_n$. These
are the elements $S\in \F_n$ such that for all $x\in [n]\setminus S$, $S\cup \{x\}\notin \F_n$.
We shall say that $\F$ is {\it partition-intersecting with $m$ components} if there exists $m$
such that for all $n\ge 2$, the set $M_n\subseteq\F_n$ of maximal elements can be partitioned into $m$
disjoint nonempty classes
$M_n = C_1\sqcup \cdots \sqcup C_m$ with the property that for all $i$, the intersection
$\bigcap_{S\in C_i} S$ is nonempty, and for all $i\ne j$ and any $S\in C_i$ and $T\in C_j$,
$S\cap T = \emptyset$.

The following lemma describes the homology groups of partition-intersecting families.

\edef\main{\the\thmcount}
\proclaim Lemma \advthm. Let $\F$ be a downward closed set family that is partition-intersecting
with $m$ components. Suppose further that $[n]\notin F_n$ for all $n\ge 2$.
Then for any cross-cut $C$ of the lattice $L = (\F\cap 2^{[n]}) \cup \top$, we have
$\tilde H_0\bigl(\Delta(C)) = \ZZ^{m-1}$ and for all $k> 0$, $\tilde H_k(\Delta(C)) = 0$.

\proof
Fix $n\ge 2$.
By hypothesis, the set $M_n$ of coatoms of $\F_n$ can be partitioned into
$M_n = C_1\sqcup C_2 \sqcup \cdots \sqcup C_m$ where each $C_i$ is nonempty, and for
all $i$, and all subsets $C'$ of $C_i$, the intersection $\bigcap_{S\in C'} S$ is nonempty. This tells
us that for any $i$, any subset of $C_i$ is not spanning and is thus a face in $\Delta(M_n)$. In other words,
the simplicial complex $\Delta(M_n)$ is the union of $m$ disjoint simplices, which has
$\tilde H_0(\Delta(M_n)) = \ZZ^{m-1}$ and $\tilde H_k(\Delta(M_n)) = 0$ for all $k > 0$.
This completes the proof.\slug

The proof of this lemma is shorter than the very definition of partition-intersecting family, and
from this one might suspect that the definition is too artificial or convoluted to be applicable.
However, we shall see
that in fact many set families arising in combinatorial number theory are partition-intersecting.
Hence Lemma~{\main} can be used to find the homology groups of the resulting lattices and in turn the homology
groups give alternating-sum identities for the cardinalities of the sets in the family.

\advsect The M\"obius function and alternating sums
\hldest{xyz}{}{proof}

This short section relates Lemma~{\main} to the cardinalities of the sets in a partition-intersecting
family.
Recall that the {\it M\"obius function} $\mu_X$ of a poset $X$ is defined on pairs of elements $(x,y)\in X^2$ with
$x\le y$, and is given by the recursive formula $\mu_X(x,x)=1$
for all $x\in X$ and
$$\mu_X(x,y) = -\sum_{x\le z<y} \mu_X(x,z).\adveq$$
The following
theorem, which appears in Rota's 1964 paper~\bref{14},
relates the M\"obius function to the homological data of a poset.

\proclaim Theorem A. Let $L$ be a lattice with bottom element $\bottom$ and
top element $\top$. For any cross-cut $C$ of $L$,
$$\mu_L(\bottom, \top) = \sum_{k=0}^n (-1)^k \rank \tilde H_k(\Delta(C)),\adveq$$
where $\tilde H_k(\Delta(C),\ZZ)$ is the $k$th reduced homology group of $\Delta(C)$.\slug

The M\"obius function generalizes the counting method
used to derive the inclusion-exclusion formula; indeed, letting
$B_n$ denote the lattice of all subsets of $[n]$, ordered by inclusion, one checks that
$\mu_{B_n}(\emptyset, S) = (-1)^{|S|}$ for any $S\in B_n$. This fact will be used in the proof of
the following proposition, which is a corollary of Lemma~{\main}. Combined with the proof of the earlier lemma,
the overall method employed is similar to
one used by M.~K.~Goh, J.~Hamdan, and J.~Saks~\bref{9} to show
that the M\"obius function of the poset of all arithmetic progressions contained in $[n]$ equals the
number-theoretic M\"obius function evaluated at $n-1$.

\edef\altsum{\the\thmcount}
\proclaim Theorem \advthm. Let $\F$ be a family of integer sets that is partition-intersecting with $m$ components.
Suppose furthermore that $\F$ is downward closed and that
for all $n\ge 2$, $\F_n = \F\cap 2^{[n]}$ does not contain $[n]$ itself.
Then, letting $F_{n,k}$ denote the number of elements
of cardinality $k$ in $\F_n$ for $0\le k<n$, we have
$$\sum_{k=1}^{n-1} (-1)^k F_{n,k} = 1-m\adveq$$
for all $n\ge 2$.

\proof Fix $n\ge 2$ and let $L_n$ be the lattice obtained by adjoining a top element to $\F_n$.
Since $\F_n$ is downward closed, for all $x\in \F_n$ the induced subposet
$$\{z\in L_n : 0\le z\le x\}$$
is isomorphic to the Boolean lattice $B_{|x|}$ and we have $\mu_{L_n}(0,x) = (-1)^{|x|}$. Then since
$[n]\notin F_n$, we have
$$\sum_{k=1}^{n-1} (-1)^k F_{n,k} = \sum_{x\in L_n\setminus\{\top\}} \mu_{L_n}(\emptyset,x)
= -\mu_{L_n}(\emptyset, \top).\adveq$$
But letting $C$ be the cross-cut of coatoms in $L_n$, by Theorem~A we have
$$\mu_{L_n}(\emptyset,\top) = \sum_{k=0}^n (-1)^k \rank \tilde H_k (\Delta(C)),$$
and since $\F$ is partition-intersecting, Lemma~{\main} tells us that
$$\sum_{k=0}^n (-1)^k \rank \tilde H_k (\Delta(C)) = \rank\tilde H_0(\Delta(C)) = m-1,$$
which is what we wanted.\slug

\advsect Primitive, pairwise coprime, and product-free sets
\hldest{xyz}{}{prim}

We now show that many families of integer sets that arise in combinatorial number theory are
partition-intersecting.
As motivation, observe that the set
$\{2,3,5,7,11,\ldots\}$ of primes possesses many pleasant properties that may be generalized
to other sets of integers. For instance, a set $S$ of integers is said to be
\medskip
\item{i)} {\it primitive} if for any two distinct $i,j\in S$, $i$ does not divide $j$;
\smallskip
\item{ii)} {\it pairwise coprime} if any two distinct $i,j\in S$ have $\gcd(i,j) = 1$; and
\smallskip
\item{iii)} {\it product-free} if for any $i,j\in S$ not necessarily distinct, $ij\notin S$.
\medskip
It is easily seen
that the set of prime numbers satisfies the requirements for all three of these definitions, and of
course, so does any finite subset of the primes. Denote $\{1,2,\ldots,n\}$ by $[n]$
for short and let $P_n$ be the number of primitive subsets of $[n]$. Likewise, let $Q_n$ be the number
of pairwise coprime subsets of $[n]$ and let $R_n$ be the number of product-free subsets of $[n]$.
These quantities have been the subject of interest of several papers over the past thirty years.
All of them were treated by a 1990 paper
of P.~J.~Cameron and P.~Erd\H{o}s~\bref{6}. They showed that
\edef\pnbounds{\the\eqcount}
$${c_1}^n \le P_n \le {c_2}^n\adveq$$
where $c_1 = 1.44967\ldots$ and $c_2 = 1.59\ldots$,
\edef\qnbounds{\the\eqcount}
$$2^{\pi(n)}e^{(1/2+o(1))\sqrt n} \le Q_n \le 2^{\pi(n)}e^{(2+o(1))\sqrt n}\adveq$$
where $\pi$ is the prime-counting function.
This result was subsequently improved to
$$ Q(n) = 2^{\pi(n)} e^{\sqrt n (1+O(\log\log n/\log n))}\adveq$$
by N.~J.~Calkin and A.~Granville~\bref{5}.
The paper of Cameron and Erd\H{o}s also gives the lower bound $R_n \ge 2^{n-\sqrt n}$,
and the authors note that any product-free sequence has an upper density of less than $1$, but there are
product-free sequences with density more than $1-\eps$ for any $\eps>0$.

Regarding primitive sets,
Cameron and Erd\H{o}s conjectured that the limit $\lim_{n\to\infty} P(n)^{1/n}$ exists; this was proven by
R.~Angelo in 2017~\bref{2}.
A separate proof was presented in a paper by H.~Liu, P.~P.~Pach, and R.~Palincza~\bref{12};
it yielded constants $c_1 = 1.571068$ and $c_2 = 1.574445$. Shortly afterward, an independent approach
of N.~McNew~\bref{13} further improved the lower bound $c_1$ to $1.572939$.

A great deal of attention has
been given to infinite primitive sets. In 1993, P.~Erd\H{o}s and Z.~Zhang~\bref{8}
proved that for any primitive set $A$,
$$\sum_{a\in A} {1\over a\log a} \le 1.84.\adveq$$
The sum $\sum_p 1/(p\log p)$ equals $1.6366\ldots$ when $p$ ranges over all primes, and Erd\H{o}s conjectured
that the bound he and Zhang gave could be improved to this value. This was proven very recently,
in a 2022 preprint of J.~D.~Lichtman~\bref{11}.

Let $P_{n,k}$ be the number of primitive subsets of $[n]$ of cardinality exactly $k$, so that
$P_n = \sum_{k=0}^n P_{n,k}$. Define $Q_{n,k}$ and $R_{n,k}$ analogously for pairwise coprime sets
and product-free sets. None of these quantities have
any known formula, as far as the authors are aware. In this section we shall prove alternating-sum
identities for these three counts. First we list some simple facts about these quantities in special cases.

\proclaim Proposition~\advthm. For $n\ge 2$ we have
\medskip
\item{i)} $P_{n,1} = Q_{n,1} = n$ and $R_{n,1} = n-1$;
\smallskip
\item{ii)} $P_{n,2} = \sum_{i=2}^n \bigl(i-d(i)\bigr)$, where $d(i)$ counts the number of divisors of $i$;
\smallskip
\item{iii)} $Q_{n,2} = \sum_{i=2}^n \varphi(i)$, where $\varphi$ is Euler's totient function;
\smallskip
\item{iv)} $R_{n,2} = {n\choose 2} - n - \lfloor \sqrt n\rfloor + 2$; and
\smallskip
\item{v)} for any $p$ prime and $k\ge 3$, $F(p,k) = F(p-1,k) + F(p-1,k-1)$, where $F$ can be any of $P$,
$Q$, or $R$.

\proof Assertion (i) is obvious, as are assertions (ii) and (iii) once one notes that in both
cases the sum runs
over the larger element in each $2$-element set. To prove (iv), we start with all $2$-element
subsets of $[n]$, then remove the $n-1$ subsets containing $1$ as well as the $\lfloor\sqrt n\rfloor -1$
subsets that contain a square and its square root.

For the final assertion, note that every primitive $k$-subset of $[p]$ either does not contain $p$ or it does.
In the first case, it is counted by $P_{p-1,k}$, and in the second case, it is counted by $P_{p-1,k-1}$,
since we can add $p$ to any primitive $(k-1)$-subset of $[p-1]$ without violating primitivity (for this
we need $k\ge 3$, otherwise the $(k-1)$-subset could be $\{1\}$). The same argument works
in the case of relatively coprime and product-free sets, since $p$ is coprime to all smaller integers.\slug

\proclaim Theorem~\advthm. The families of primitive sets, pairwise coprime sets, and product-free sets
are all partition-intersecting, with two components in the first case and
one component in the other two cases.

\proof
Fix $n\ge 2$ and let $p$ be the largest prime less than or equal to $n$. By the
famous 1852 theorem of P.~Chebyshev~\bref{7},
we have $2p>n$. This will be important in two of the cases. We now show that each of these
families is partition-intersecting, computing the number $m$ of components along the way.

We deal with primitive sets first. In this case, the
collection $C$ of maximal primitive sets contains the singleton $\{1\}$, but no other element of $C$
contains $1$, since $1$ divides every other positive integer. So we set $C_1 = \{\{1\}\}$
and $C_2 = C\setminus \{1\}$. To any subset of $\{2,3,\ldots,n\}$, we may add the element $p$
without violating primitivity, since no element divides $p$ and $2p>n$. Thus every element of $C_2$
contains $p$ and we see that $m = 2$.

For the family of pairwise coprime sets, every maximal set contains $1$, since $\gcd(1,i) = 1$
for all $1\le i\le n$. So we have $C_1 = C$ and $m = 1$.

Lastly, we note that no product-free set contains $1$, since $1\cdot 1 = 1$,
and every maximal product-free subset of $[n]$ must also contain $p$, since $p$ is not the product
of any two elements of $\{2,3,\ldots,n\}$ and $2p>n$. Hence again we have $C_1 = C$ and $m=1$
in this case as well.\slug

The interval $[n]$ is neither primitive, nor pairwise coprime (for $n\ge 4$),
nor product-free. Furthermore,
these three properties are all closed under taking subsets, so the families are downward closed,
which gives us the following corollary of Theorem~{\altsum}. (Note that since $[2]$ and
$[3]$ are pairwise coprime, the cases $n=2$ and $n=3$ must
be verified separately in the case of pairwise coprime sets.)

\proclaim Corollary~\advthm. For $n\ge 2$, we have
\medskip
\item{i)} $\sum_{k=0}^n (-1)^k P_{n,k} = -1$;
\smallskip
\item{ii)} $\sum_{k=0}^n (-1)^k Q_{n,k} = 0$; and
\smallskip
\item{iii)} $\sum_{k=0}^n (-1)^k R_{n,k} = 0$.\noskipslug
\medskip

Values of $P_{n,k}$, $Q_{n,k}$, and $R_{n,k}$ for small $n$ and $k$ are presented in the appendix.
For the purpose of exposition, we have chosen to focus on three specific set families in this section,
but many more multiplicative
conditions can be placed on integer sets to yield partition-intersecting families. Two more examples
from~\bref{6} are treated by the following proposition.

\edef\moreexamples{\the\thmcount}
\proclaim Proposition \advthm. If $\F$ is the family of all finite subsets $S$ of $\NN$ such that
\medskip
\item{i)} for all $i,j,k,l\in S$ all distinct, $ij\ne kl$;
\smallskip
\item{ii)} for all $i,j,k\in S$ with $i\notin \{j,k\}$, $i$ does not divide $jk$; or
\smallskip
\item{iii)} for all $i,j\in S$ with $i<j$, $i$ divides $j$,
\medskip\noindent
then $\F$ is partition-intersecting with one component and, letting $F_{n,k}$ denote the
number of sets in $\F\cap 2^{[n]}$ of cardinality $k$, we have
$$\sum_{k=0}^{n} (-1)^k F_{n,k} = 0,$$
for $n$ large enough depending on the case.

\proof Cases (i) and (ii) are proved in a fashion similar to the case of product-free sets. Case (iii)
is similar to the case of pairwise coprime sets. The details are left to the reader.\slug

Note that if $m$ is a constant not depending on $n$, then for a family to be partition-intersecting
with $m$ components we must have $m$ equal to either $1$ or $2$, since there are at most two maximal
elements in $\F_2$. But by relaxing $n\ge 2$ to $n\ge N$ for some larger $N$ in the definition of
partition-intersecting, the value of $m$ would then range between $1$ and ${N\choose N/2}$. One could also
imagine extending the definition to let $m$ vary as a function of $n$, which may allow more set
families to be treated.

Lastly, note that for all of these examples,
the alternating-sum identities themselves are not so difficult to prove from
scratch. For six of the seven families mentioned above,
let $p$ be a prime in the interval $(n/2,n]$ and note that any $S\subseteq [n]\setminus \{p\}$ is
in the family if and only if $S\cup \{p\}$ is. The one exception is in the case of primitive sets,
when $S = \{1\}$ is primitive but $\{1,p\}$ is not, which explains why the alternating sum equals $-1$ in this
case.

\advsect Coprime-free sets
\hldest{xyz}{}{coprimefree}

In part (iii) of Proposition~{\moreexamples} above, we negated the condition ``$i$ does not divide $j$''
from the definition of primitive sets and ended up with another partition-intersecting family. In this
section we show what happens if one negates the condition $\gcd(i,j) = 1$, that is, if we require
that $\gcd(i,j) > 1$ for all $i\ne j$ in each of our sets. These sets are said to be {\it coprime-free},
and they are not partition-intersecting in general.
The asymptotic number of coprime-free sets was
treated by the paper~\bref{6} of Cameron and Erd\H{o}s, and subsequently in the
paper~\bref{5}
by Calkin and Granville. (Both papers were mentioned in the previous section.)
The latter paper showed that the number
of subsets $S\subseteq [n]$ for which $\gcd(i,j) \ne 1$ for all distinct $i,j\in S$ is
$$ 2^{\lfloor n/2\rfloor} + 2^{\lfloor n/2\rfloor-N} + O\biggl( 2^{\lfloor n/2\rfloor - N}
\exp \biggl( -C {n\over \log^2 n\log\log n}\biggr)\biggr),\adveq$$
for an absolute constant $C>0$, where $N = \bigl(e^{-\gamma}+o(1)\bigr) n/\log\log n$.
(Here $\gamma$ denotes the Euler-Mascheroni constant $0.5772\ldots$.)

In the proof of the proposition
below we will make use of the notion of a nerve complex, which we define here, alongside some
more terminology concerning simplicial complexes.
Let $C = \{U_i\}_{i\in I}$
be a set family. The {\it nerve complex} of $C$ is a simplical complex defined on the vertex set $I$,
where a set $J\subseteq I$ forms a face if and only if $\bigcap_{i\in J} U_i \ne \emptyset$. Borsuk's nerve
theorem~\bref{4} states that for
a simplicial complex $\Delta$ and a set $\{\Delta_i\}_{i\in I}$
of subcomplexes of $\Delta$, if
\medskip
\item{i)} $\bigcup_{i\in I}\Delta_i = \Delta$; and
\smallskip
\item{ii)} any nonempty intersection of finitely many $\Delta_i$ is contractible,
\medskip\noindent
then the nerve complex of $\{\Delta_i\}_{i\in I}$ is
homotopy equivalent to $\Delta$.

The {\it simplicial star} of a vertex $x$ in a simplicial complex is the set of all faces
containing $x$.
The {\it link} of $x$ is the set of all faces $F$ such that $x\notin F$ but $F\cup \{x\}$
is a face.

Also of use to us is the
the Mayer--Vietoris sequence, which for reduced simplicial homology says that
$$\eqalign{
\cdots \to \tilde H_s(A)\oplus \tilde H_s(B) \to \tilde H_s(X)\to\tilde H_{s-1}(A\cap B) \to
\tilde H_{s-1}(A)\oplus \tilde H_{s-1}(B)\to\cr
\cdots\to \tilde H_0(A\cap B) \to \tilde H_0(A)\oplus \tilde H_0(B) \to \tilde H_0(X) \to 0,\cr
}\adveq$$
whenever $A$ and $B$ are simplicial subcomplexes whose union is
the simplicial complex $X$ (see, e.g., Remark~5.18 of \bref{10}).

The following proposition computes the zeroth homology of the posets corresponding to coprime-free subsets
of $[n]$, and shows that the first homology is trivial.

\edef\coprimefree{\the\thmcount}
\proclaim Proposition~\advthm. Let $\F$ be the family of finite coprime-free subsets of $\NN$. Let $\F_n$
denote $\F\cap 2^{[n]}$ and let $\Delta_n$ be the simplicial complex corresponding to the
cross-cut of coatoms in the lattice $\F_n\cup \top$.
Then we have
$$\tilde H_0(\Delta_1) = 0,\qquad \tilde H_0(\Delta_2) = \ZZ, \qquad\tilde H_0(\Delta_3) = \ZZ^2,\adveq$$
and for $n\ge 4$ we have $\tilde H_0(\Delta_n) = \ZZ^{C(n)+1}$, where $C(n)$ is the number of primes
in the range $(n/2,n]$. Furthermore, $\tilde H_1(\Delta_n)$ is trivial.

\proof The assertions for $1\le n\le 3$ easily follow from the observation that the only nonempty coprime-free
subsets of $\{1,2,3\}$ are the singleton sets $\{1\}$, $\{2\}$, and $\{3\}$.

Now fix $n\ge 4$. Any prime $p$ in the range $(n/2,n]$ is coprime to every integer in $[n]\setminus\{p\}$,
so the set $\{p\}$ is a maximal coprime-free subset of $[n]$ and hence is an isolated vertex of $\Delta_n$.
The set $\{1\}$ is also an isolated vertex for the same reason, giving us $C(n)+1$ isolated points in $\Delta_n$.
To prove that $\tilde H_0(\Delta_n) = \ZZ^{C(n)+1}$ we must show that the rest of $\Delta_n$ forms one
nonempty connected
component. Note first that for every prime $p\in [2,n/2]$ (there are such primes because $n\ge 4$),
the set $x_p = p\NN\cap [n]$ is a non-singleton
maximal coprime-free subset of $[n]$. Since $2p\in x_p$ for all $p\in (2,n/2]$, $x_p\cap x_2 \ne \emptyset$
for all these coatoms, and thus these vertices are connected by some nontrivial face in $\Delta_n$.
Lastly, consider an arbitrary vertex of $\Delta_n$ that is not one of the $C(n)+1$ singletons described
above. This coatom must contain some multiple of a prime $p$ in $[2,n/2]$. Hence $x\cap x_p\ne \emptyset$.

We will now prove that $\tilde H_1(\Delta_n)$ is trivial for all $n$. To do so, we will
use Borsuk's nerve theorem to reduce the problem to studying a simpler simplicial complex.
For all $i\in [n]$, let $U_i$ be the
simplicial subcomplex of $\Delta_n$ induced on the vertices that contain $i$.
Note that the set $\{U_i\}_{i\in [n]}$ covers $\Delta_n$. This is true because faces of $\Delta_n$ are present
wherever a set of maximal coprime-free subsets have nonempty intersection, and this intersection must contain
some element of $[n]$. Moreover, for any $J\subseteq [n]$, let
$U_J = \bigcap_{i\in J} U_i$. A face of $\Delta_n$ belongs to this intersection precisely when
each of its vertices contains every element of $J$, so this is
the subcomplex of $\Delta_n$ induced on vertices $x$ such
that $J\subseteq x$. If $J$ is not coprime-free,
then clearly $U_J$ is empty, and if $J$ is coprime-free and nonempty,
then since all vertices in $U_J$ are supersets of $J$, they all intersect and thus the complex $U_J$ is a simplex
and {\it a fortiori} contractible.
Applying Borsuk's nerve theorem, we have shown that $\Delta_n$ is homotopy equivalent to the simplicial
complex $X_n$ whose vertex set is $[n]$ and in which a subset $J\subseteq [n]$ is a face if and only if
$J$ is coprime-free.

We proceed by induction on $n$. That $\tilde H_1(X_n) = 0$ for $1\le n\le 4$
can easily be verified by hand. Now let
$n>4$. If $n$ is prime, then $X_n$ is simply the union of $X_{n-1}$ and a new isolated vertex,
so $\tilde H_0(X_n) = \tilde H_0(X_{n-1}) \oplus \ZZ$
but all other homologies remain unchanged, and by the induction hypothesis,
$\tilde H_1(X_n) = 0$.

The interesting case is when $n>4$ is not prime. We let $A$ be the simplicial star of the vertex
$n$ in $X_n$, and let $B = X_{n-1}$, considered as a subset of $X_n$.
It is clear that the union of $A$ and $B$ is $X_n$,
since $A$ contains all faces containing $n$ and $B$ contains all other faces.
The relevant section
of the Mayer--Vietoris sequence is
$$\eqalign{
\cdots &\to \tilde H_1(A)\oplus \tilde H_1(B) \to \tilde H_1(X_n) \to \tilde H_0(A\cap B) \to\cr
&\qquad\qquad\qquad\qquad\qquad \tilde H_0(A)\oplus  \tilde H_0(B) \to  \tilde H_0(X_n) \to 0 \cr
}$$
Note that $A$, being a simplicial star, is contractible. Moreover, we proved above that
$\tilde H_0(X_n) = C(n)+1$ and $\tilde H_0(B) = C(n-1)+1$,
and by the induction hypothesis, $\tilde H_1(B) = 0$.

Next we count the connected components in $A\cap B$, whose vertices are the integers in $[n-1]$ that share
a prime factor with $n$.
Let $i,j\in A\cap B$, let $p$ be a prime that divides both $i$ and $n$, and let $q$ be a prime that
divides both $j$ and $n$. Writing $n = kp = lq$, we have $k,l\ge 2$, since $n$ is not prime. If, furthermore,
$n/2$ is not prime, then $k,l\ge 3$, so both $2p$ and $2q$ are in $A\cap B$ and we have the path
$i \to 2p \to 2q \to j$ along the $1$-skeleton of $A\cap B$. In the case that $n/2$ is prime, the case
that $k$ or $l$ is equal to $2$ corresponds to either $i$ or $j$ being $n/2$, and we already
know that the prime $n/2$ is an isolated vertex of $X_{n-1}$, so it is
also an isolated vertex of $A\cap B$.
We have shown that
$$\rank\tilde H_0(A\cap B) = \one_{n/2\hbox{\sevenrm\ is prime}}.$$

Putting all these facts into the exact sequence above, we have
$$\cdots \to 0 \to \tilde H_1(X_n)
\to \ZZ^{\one_{n/2\hbox{\fiverm\ is prime}}} \to 0 \oplus  \ZZ^{C(n-1)+1}\to  \ZZ^{C(n)+1} \to 0$$
But since $n$ is composite, we have
$$C(n-1) = \cases{C(n)+1, & if $n/2$ is prime;\cr C(n), & otherwise.}$$
Thus we see that $C(n-1) + 1 = C(n)+1 + \one_{n/2\hbox{\sevenrm\ is prime}}$, and can
conclude that $\rank\tilde H_1(X_n) = 0$.\slug

This proof illustrates that the situation with coprime-free sets is rather different from the partition-intersecting
families we considered before. From Proposition~{\coprimefree}
we cannot conclude anything about the alternating sum corresponding to coprime-free sets, and next we will
show that the proposition cannot be improved in general.

In the following description of $\Delta_{143}$, which was communicated to us by M.~Adamaszek~\bref{1},
we shall use terminology established in the work of J.~A.~Barmak and
E.~G.~Minian~\bref{3}. A vertex $x$ in a simplicial complex $X$ is said to be {\it dominated}
by another vertex $y$ if the link of $x$ is the simplicial star of $y$ in the subcomplex
$X\setminus \{x\}$. In other words, for every face $F$ of $X$ that contains $x$, $F\cup \{y\}$ is also a face
of $X$. Hence removing the vertex $x$
from the complex does not change the homotopy type, and the replacement of $X$ with $X\setminus\{x\}$
is called an {\it elementary strong collapse} of $X$.

\edef\counterexample{\the\thmcount}
\proclaim Proposition~\advthm. Let $\F$ be the family of finite coprime-free subsets of $\NN$. Let $\F_{143}$
denote $\F\cap 2^{[143]}$ and let $\Delta_{143}$ be the simplicial complex corresponding to the
cross-cut of coatoms in the lattice $\F_{143}\cup \top$. We have $\tilde H_2(\Delta_{143}) \ne 0$.

\proof By Borsuk's nerve theorem as we employed it in the proof of Proposition~{\coprimefree},
$\Delta_{143}$ is homotopy equivalent to the simplicial complex $X_{143}$ whose vertex set is $[143]$
and in which a nonempty subset $S\subseteq [143]$ forms a face if and only if $S$ is coprime-free.

We now describe a sequence of elementary strong
collapses. We check every pair $(i,j)\in [143]$ with $i\ne j$ to see if $j$ dominates $i$. For each such pair we may
remove $i$ from the complex without altering the homology groups.
Note that in our context, $j$ dominates $i$ if and only
if the set of primes dividing $i$ is a subset of the set of primes dividing $j$.
It is possible that both $i$ dominates $j$ and $j$ dominates $i$ (this happens if they have the same
set of distinct prime factors with different multiplicities). In this case it is valid to remove either $i$ or $j$,
but we shall make the choice to eliminate the larger integer.
When the algorithm
materminates, we are left with a simplicial complex $X'$ whose vertex set is the set of all maximal squarefree
integers in $[143]$, that is, squarefree integers $i$ such that no squarefree multiple of $i$ belongs to $[143]$.

We can safely ignore vertices in $X'$ corresponding to primes in $(143/2, 143]$ as these isolated vertices
have no bearing on $\tilde H_2(X')$. The main component in $X'$
is the simplicial complex on the vertex set
$$\eqalign{
\{&30, 42, 58, 62, 66, 70, 74, 77, 78, 82, 85, 86, 87, 91, 93, \cr
&94, 95, 102, 105, 106, 110, 111, 114, 115, 118, 119, \cr
&122, 123, 129, 130, 133, 134, 138, 141, 142, 143 \} \cr
}$$
where, again, a set of vertices forms a face if and only if it is coprime-free.
Consider the subcomplex on the vertices $\{42, 66, 77, 78, 91, 143\}$. From the factorizations
$$\eqalign{
42 &= 2\cdot 3\cdot 7,\qquad\hskip 4.8pt 66 = 2\cdot 3\cdot 11, \qquad\hskip28.8pt 77 = 7\cdot 11,\cr
78 &= 2\cdot 3\cdot 13, \qquad 91=7\cdot 13, \qquad\hbox{and}\qquad 143 = 11\cdot 13, \cr
}$$
it is easy to check by hand that these vertices form the surface of an octahedron. For this to be a nontrivial
generator of $\tilde H_2(X')$, we must rule out the possibility that
it is the boundary of a larger cycle. The complex $X'$ is small enough that this
is amenable to machine-assisted verification.
\slug

In Appendix B, we include a Sage program that performs the computation alluded to in the proof of Proposition~8.
This program reports that $\tilde H_2(X_{143}) = \ZZ$, and $\tilde H_2(X_n) = 0$ for all $n\le 142$.

\advsect A generalisation of primitive sets
\hldest{xyz}{}{general}

A primitive set is a set $S$ that does not contain more than one multiple of any integer in the set.
We now generalize the condition to forbid having more than $s$ multiples of an integer in the set. We
will refer to sets satisfying this condition as {\it $s$-multiple sets}, and when
$s = 1$ we recover the definition of a primitive set. In the more general case, we still have the
useful property that any maximal $s$-multiple set contains $1$ or it contains
every prime in the range $(n/2,n]$.
The following theorem describes the homology of lattices associated with $s$-multiple
set families.

As with coprime-free sets, $s$-multiple families are not partition-intersecting in general, but
unlike in the previous section, we are
able to fully describe the homology of the simplicial complexes corresponding to $s$-multiple set families.

\edef\smultiple{\the\thmcount}
\proclaim Theorem~\advthm. Let $s\ge 2$ be an integer and let $\F_s$ be the family
of $s$-multiple subsets of $\NN$. Fix $n>s$, let $\F_{s,n} = \F_s \cap 2^{[n]}$, and let
$\Delta_n$ be the simplicial complex formed by the cross-cut of coatoms in $L_{s,n} = \F_{s,n}\cup \top$.
Then we have
$$\rank\tilde H_t (\Delta_n) = \cases{ {n-2\choose s-1},& if $t = s-1$;\cr 0,& otherwise.\cr}\adveq$$

\proof Borsuk's nerve theorem once again tells us that it suffices to study the simplicial complex
$X_n$ on the vertex set $\{1,\ldots,n\}$, with a face for every $s$-multiple set of vertices.
We again use the Mayer--Vietoris sequence, so let
$A$ be the simplicial star of $1$ in $X_n$, and let $B$ be the full simplicial subcomplex on
$\{2,\ldots,n\}$. The intersection $A\cap B$ is the link of $1$, in which every subset of
$\{2,\ldots,n\}$ of cardinality $s-1$ is a maximal face,
since $1$ may be added to any such set without violating the
$s$-multiple condition, but any subset of $\{2,\ldots,n\}$ of cardinality $s$ already contains $s$
multiples of $1$. Thus $A\cap B$ is the $(s-2)$-skeleton of an $(n-2)$-simplex, and we have
$$\rank\tilde H_t (A\cap B) = \cases{ {n-2\choose s-1},& if $t = s-2$;\cr 0,& otherwise.\cr}\adveq$$

The subcomplex $A$, being the simplicial star of $1$, is evidently contractible.
Now let $p$ be a prime in the range $(n/2,n]$. Since $1$ is not a vertex of $B$, and since
any multiple of $p$ is greater than $n$, we see that for any $F$ in $B$ such that
$p\ne F$, the set $F\cup \{p\}$ is also a face of $B$. This means that the entire simplicial
complex $B$ is the simplicial star of $p$, and is thus contractible as well. But since
$\tilde H_k(A)\oplus \tilde H_k(B) \cong 0$ for all $k$, the Mayer--Vietoris sequence informs us that
$$\tilde H_k(X_n) \cong \tilde H_{k-1}(A\cap B)\adveq$$
for all $k\ge 1$. This completes the proof.\slug

Finally, by reasoning as in the proof of Theorem~{\altsum}, we have a corollary concerning alternating sums of
counts of $s$-primitive sets.

\proclaim Corollary~\advthm. Let $s\ge 2$ be an integer, let $n>s$,
and let $P_{s,n,k}$ denote the number of $s$-multiple subsets of $[n]$ with
cardinality exactly $k$. Then
$$ \sum_{k=0}^n (-1)^k P_{s,n,k} = (-1)^s{n-2\choose s-1} .\noskipslug\adveq$$

\section Acknowledgements
\hldest{xyz}{}{acks}

The first author is grateful for the generous financial support of Luc Devroye, and the second author is
funded by the Natural Sciences and Engineering Research Council of Canada. We thank the users of MathOverflow
for their help in correcting an error that appeared in a preliminary version of the paper;
particular thanks are due to user Andy Putman who suggested the use of the nerve complex,
and to user Micha{\l} Adamaszek for communicating the example of a coprime-free set family
with nontrivial second homology. We also thank Andrew Granville, Jacob Reznikov,
and Daniel Wise for helpful technical suggestions and enjoyable discussions.
Last but certainly not least, we are indebted to the anonymous referee for giving our work a careful read,
and for finding the aforementioned error in an early version of our paper.

\section References
\hldest{xyz}{}{refs}

\parskip=0pt
\hyphenpenalty=-1000 \pretolerance=-1 \tolerance=1000
\doublehyphendemerits=-100000 \finalhyphendemerits=-100000
\frenchspacing
\def\bref#1{[#1]}
\def\beginref{\noindent}
\def\endref{\medskip}
\vskip\parskip

\beginref
\parindent=20pt\item{\bref{1}}
\hldest{xyz}{}{bib1}%
Micha{\l} Adamaszek,
``Amending flawed `proof' that homology groups are zero,'' {\sl URL (version: 2023-09-24):} {\tt https://mathoverflow.net/q/455210}.
\endref
\beginref
\parindent=20pt\item{\bref{2}}
\hldest{xyz}{}{bib2}%
Rodrigo Angelo,
``A Cameron and Erd\H{o}s conjecture on counting primitive sets,''
{\sl Integers: Electronic Journal of Combinatorial Number Theory}\/
{\bf 18}
(2018),
article no.~A25.
\endref
\beginref
\parindent=20pt\item{\bref{3}}
\hldest{xyz}{}{bib3}%
Jonathan Ariel Barmak
and Elias Gabriel Minian,
``Strong homotopy types, nerves, and collapses,''
{\sl Discrete and Computational Geometry}\/
{\bf 47}
(2012),
301--328.
\endref
\beginref
\parindent=20pt\item{\bref{4}}
\hldest{xyz}{}{bib4}%
Karol Borsuk,
``On the imbedding of systems of compacta in simplicial complexes,''
{\sl Fundamenta Mathematicae}\/
{\bf 35}
(1948),
217---234.
\endref
\beginref
\parindent=20pt\item{\bref{5}}
\hldest{xyz}{}{bib5}%
Neil James Calkin
and Andrew Granville,
``On the number of co-prime-free sets,'' {\sl Number Theory} (New York, 1991--1995), 9--18.
\endref
\beginref
\parindent=20pt\item{\bref{6}}
\hldest{xyz}{}{bib6}%
Peter Jephson Cameron
and Paul Erd\H{o}s,
``On the number of sets of integers with various properties,'' {\sl Number Theory} (Banff, 1990), 61--79.
\endref
\beginref
\parindent=20pt\item{\bref{7}}
\hldest{xyz}{}{bib7}%
Pafnuty Chebyshev,
``M\'emoire sur les nombres premiers,''
{\sl Journal de math\'ematiques pures et appliqu\'ees, S\'erie 1}\/
{\bf 17}
(1852),
366--390.
\endref
\beginref
\parindent=20pt\item{\bref{8}}
\hldest{xyz}{}{bib8}%
Paul Erd\H{o}s
and Zhenxiang Zhang,
``Upper bound of $\sum 1/(a_i \log a_i)$ for primitive sequences,''
{\sl Proceedings of the American Mathematical Society}\/
{\bf 117}
(1993),
891--895.
\endref
\beginref
\parindent=20pt\item{\bref{9}}
\hldest{xyz}{}{bib9}%
Marcel Kieren Goh,
Jad Hamdan,
and Jonah Saks,
``The lattice of arithmetic progressions,''
{\sl Australasian Journal of Combinatorics}\/
{\bf 84}
(2022),
357--374.
\endref
\beginref
\parindent=20pt\item{\bref{10}}
\hldest{xyz}{}{bib10}%
Dmitry Kozlov,
{\sl Combinatorial Algebraic Topology}
(Berlin:
Springer-Verlag,
2008).
\endref
\beginref
\parindent=20pt\item{\bref{11}}
\hldest{xyz}{}{bib11}%
Jared Duker Lichtman,
``A proof of the Erd\H{o}s primitive set conjecture,''
{\sl arXiv:2202.02384}\/
(2022),
22 pp.
\endref
\beginref
\parindent=20pt\item{\bref{12}}
\hldest{xyz}{}{bib12}%
Hong Liu,
P\'eter P\'al Pach,
and Rich\'ard Palincza,
``The number of maximum primitive sets of integers,''
{\sl Combinatorics, Probability and Computing}\/
{\bf 30}
(2021),
781--795.
\endref
\beginref
\parindent=20pt\item{\bref{13}}
\hldest{xyz}{}{bib13}%
Nathan McNew,
``Counting primitive subsets and other statistics of the divisor graph of $\{1,2,\ldots,n\}$,''
{\sl European Journal of Combinatorics}\/
{\bf 92}
(2021),
article no.~103237.
\endref
\beginref
\parindent=20pt\item{\bref{14}}
\hldest{xyz}{}{bib14}%
Gian-Carlo Rota,
``On the foundations of combinatorial theory I. Theory of M\"obius Functions,''
{\sl Zeitschrift f\"ur Wahrscheinlichkeitstheorie und verwandte Gebiete}\/
{\bf 2}
(1964),
340--368.
\endref
\beginref
\vfill\eject

\section A. Numerical tables
\hldest{xyz}{}{numerical}

This appendix contains values of $P_{n,k}$, $Q_{n,k}$,
and $R_{n,k}$ for small values of $n$ and $k$.
The row sums of $P_{n,k}$, which we called $P_n$ in section 3, appear as A051026
in the On-line Encyclopedia of Integer Sequences ({\mc OEIS}).
The row sums of $Q_{n,k}$ and $R_{n,k}$ were referred to as $Q_n$
and $R_n$ in section 3 and
appear as A084422 and A326489, respectively.
The triangular numbers $P_{n,k}$, $Q_{n,k}$, and $R_{n,k}$ are A355145, A355146, and A355147
in the {\mc OEIS}, respectively.

\vskip 100pt

\midinsert
$$\vcenter{\vbox{
\eightpoint
\centerline{\smallheader Table 1}
\medskip
\centerline{THE NUMBER $P_{n,k}$ OF PRIMITIVE $k$-SUBSETS OF $[n]$}
}}$$
\vskip-10pt
$$\vcenter{\vbox{
\eightpoint
\hrule
\medskip
\tabskip=.7em plus.2em minus .5em
\halign{
   $\hfil#$  &  $\hfil#$ & $\hfil#$ & $\hfil#$ & $\hfil#$ & $\hfil#$ & $\hfil#$ & $\hfil#$ &
   $\hfil#$ & $\hfil#$ & $\hfil#$ & $\hfil#$ & $\hfil#$ \cr
   n   & P_{n,0} & P_{n,1} & P_{n,2} & P_{n,3} & P_{n,4} & P_{n,5} & P_{n,6} & P_{n,7} & P_{n,8} & P_{n,9}\cr
   \noalign{\medskip}
   \noalign{\hrule}
   \noalign{\medskip}
1  & 1 & 1  & 0 & 0 & 0 & 0 & 0 & 0 & 0 & 0      \cr
2  & 1 & 2  & 0 & 0 & 0 & 0 & 0 & 0 & 0 & 0      \cr
3  & 1 & 3  & 1 & 0 & 0 & 0 & 0 & 0 & 0 & 0      \cr
4  & 1 & 4  & 2 & 0 & 0 & 0 & 0 & 0 & 0 & 0      \cr
5  & 1 & 5  & 5 & 2 & 0 & 0 & 0 & 0 & 0 & 0      \cr
6  & 1 & 6  & 7 & 3 & 0 & 0 & 0 & 0 & 0 & 0      \cr
7  & 1 & 7  & 12 & 10 & 3 & 0 & 0 & 0 & 0 & 0      \cr
8  & 1 & 8  & 16 & 15 & 5 & 0 & 0 & 0  & 0 & 0     \cr
9  & 1 & 9  & 22 & 26 & 13 & 2 & 0 & 0      & 0 & 0 \cr
10 & 1 & 10 & 28 & 38 & 22 & 4 & 0 & 0      & 0 & 0 \cr
11 & 1 & 11 & 37 & 66 & 60 & 26 & 4 & 0      & 0 & 0 \cr
12 & 1 & 12 & 43 & 80 & 76 & 35 & 6 & 0      & 0 & 0 \cr
13 & 1 & 13 & 54 & 123 & 156 & 111 & 41 & 60 & 0 & 0      \cr
14 & 1 & 14 & 64 & 161 & 227 & 180 & 74 & 12 & 0 & 0      \cr
15 & 1 & 15 & 75 & 206 & 323 & 299 & 161 & 47 & 6  &0     \cr
16 & 1 & 16 & 86 & 253 & 425 & 421 & 242 & 75 & 10  &0    \cr
17 & 1 & 17 & 101 & 339 & 678 & 846 & 663 & 317 & 85 & 10 \cr
 \noalign{\medskip}
 \noalign{\hrule}
    }
}}$$
\endinsert
\vfill\eject
\midinsert
$$\vcenter{\vbox{
\eightpoint
\centerline{\smallheader Table 2}
\medskip
\centerline{THE NUMBER $Q_{n,k}$ OF PAIRWISE COPRIME $k$-SUBSETS OF $[n]$}
}}$$
\vskip-10pt
$$\vcenter{\vbox{
\eightpoint
\hrule
\medskip
\tabskip=.7em plus.2em minus .5em
\halign{
   $\hfil#$  &  $\hfil#$ & $\hfil#$ & $\hfil#$ & $\hfil#$ & $\hfil#$ & $\hfil#$ & $\hfil#$ &
   $\hfil#$ & $\hfil#$ & $\hfil#$ & $\hfil#$ & $\hfil#$ \cr
   n   & Q_{n,0} & Q_{n,1} & Q_{n,2} & Q_{n,3} & Q_{n,4} & Q_{n,5} & Q_{n,6} & Q_{n,7} & Q_{n,8} \cr
   \noalign{\medskip}
   \noalign{\hrule}
   \noalign{\medskip}
1  & 1 & 1  & 0 & 0 & 0 & 0 & 0 & 0 & 0 \cr
2  & 1 & 2  & 1 & 0 & 0 & 0 & 0 & 0 & 0 \cr
3  & 1 & 3  & 3 & 1 & 0 & 0 & 0 & 0 & 0 \cr
4  & 1 & 4  & 5 & 2 & 0 & 0 & 0 & 0 & 0 \cr
5  & 1 & 5  & 9 & 7 & 2 & 0 & 0 & 0 & 0 \cr
6  & 1 & 6  & 11 & 8 & 2 & 0 & 0 & 0 & 0 \cr
7  & 1 & 7  & 17 & 19 & 10 & 2 & 0 & 0 & 0 \cr
8  & 1 & 8  & 21 & 25 & 14 & 3 & 0 & 0 & 0 \cr
9  & 1 & 9  & 27 & 37 & 24 & 6 & 0 & 0 & 0 \cr
10 & 1 & 10 & 31 & 42 & 26 & 6 & 0 & 0 & 0 \cr
11 & 1 & 11 & 41 & 73 & 68 & 32 & 6 & 0 & 0 \cr
12 & 1 & 12 & 45 & 79 & 72 & 33 & 6 & 0 & 0 \cr
13 & 1 & 13 & 57 & 124 & 151 & 105 & 39 & 6 & 0 \cr
14 & 1 & 14 & 63 & 138 & 167 & 114 & 41 & 6 & 0 \cr
15 & 1 & 15 & 71 & 159 & 192 & 128 & 44 & 6 & 0 \cr
16 & 1 & 16 & 79 & 183 & 228 & 157 & 56 & 8 & 0 \cr
17 & 1 & 17 & 95 & 262 & 411 & 385 & 213 & 64 & 8 \cr
 \noalign{\medskip}
 \noalign{\hrule}
    }
}}$$
\endinsert


\midinsert
$$\vcenter{\vbox{
\eightpoint
\centerline{\smallheader Table 3}
\medskip
\centerline{THE NUMBER $R_{n,k}$ OF PRODUCT-FREE $k$-SUBSETS OF $[n]$}
}}$$
\vskip-10pt
$$\vcenter{\vbox{
\eightpoint
\hrule
\medskip
\tabskip=.7em plus.2em minus .5em
\halign{
   $\hfil#$  &  $\hfil#$ & $\hfil#$ & $\hfil#$ & $\hfil#$ & $\hfil#$ & $\hfil#$ & $\hfil#$ &
   $\hfil#$ & $\hfil#$ & $\hfil#$ & $\hfil#$ & $\hfil#$ \cr
   n   & R_{n,0} & R_{n,1} & R_{n,2} & R_{n,3} & R_{n,4} & R_{n,5} & R_{n,6} & R_{n,7} & R_{n,8} & R_{n,9}\cr
   \noalign{\medskip}
   \noalign{\hrule}
   \noalign{\medskip}
1  &1 & 0  &  0 & 0 & 0 & 0 & 0 & 0 & 0 & 0 \cr
2  &1 & 1  &  0 & 0 & 0 & 0 & 0 & 0 & 0 & 0 \cr
3  &1 & 2  &  1 & 0 & 0 & 0 & 0 & 0 & 0 & 0 \cr
4  &1 & 3  &  2 & 0 & 0 & 0 & 0 & 0 & 0 & 0 \cr
5  &1 & 4  &  5 & 2 & 0 & 0 & 0 & 0 & 0 & 0 \cr
6  &1 & 5  &  9 & 6 & 1 & 0 & 0 & 0 & 0 & 0 \cr
7  &1 & 6  &  14 & 15 & 7 & 1 & 0 & 0 & 0 & 0 \cr
8  &1 & 7  &  20 & 29 & 22 & 8 & 1 & 0 & 0 & 0 & \cr
9  &1 & 8  &  26 & 43 & 38 & 17 & 3 & 0 & 0 & 0 & \cr
10 &1 & 9  &  34 & 68 & 76 & 47 & 15 & 2 & 0 & 0 & \cr
11 &1 & 10 &  43 & 102 & 144 & 123 & 62 & 17 & 2 & 0  \cr
12 &1 & 11 &  53 & 143 & 234 & 238 & 149 & 55 & 11 & 1 \cr
 \noalign{\medskip}
 \noalign{\hrule}
    }
}}$$
\endinsert

\vfill\eject

\def\_{\char'137}%
\section B. Computer verification of Proposition 8
\hldest{xyz}{}{verification}

The following Sage program verifies that $n=143$ is the first example for which $X_n$, as defined in Proposition~8,
has nontrivial second homology, and in fact it shows that $\tilde H_2(X_{143}) = \ZZ$. We did not pay particular
attention to the efficiency of runtime; nonetheless the entire computation took less than
eight minutes on a standard laptop computer with 16~{\mc GB} of {\mc RAM} and a 2.4~GHz Intel Core~i5 processor.

\begincode
def coprime\_free\_second\_homology(n):
\tab F = set()
\tab {\sl\# build set of maximal squarefree integers in $[n]$}
\tab for i in reversed(range(1,n+1)):
\tab \tab if not is\_squarefree(i):
\tab \tab \tab continue
\tab \tab dominated = False
\tab \tab for j in F:
\tab \tab \tab if j \% i == 0:
\tab \tab \tab \tab dominated = True
\tab \tab \tab \tab break
\tab \tab if not dominated:
\tab \tab \tab F.add(i)
\tab {\sl\# build graph on composite numbers in $F$,}
\tab {\sl\# where $i \mathbin{-\!\!\!-} j$ if and only if $\gcd(i,j) > 1$}
\tab G = Graph()
\tab for i in F:
\tab \tab if not is\_prime(i):
\tab \tab \tab G.add\_vertex(i)
\tab for i in G.vertices():
\tab \tab for j in G.vertices():
\tab \tab \tab if i < j and gcd(i,j) > 1:
\tab \tab \tab \tab G.add\_edge(i,j)
\tab S = G.clique\_complex()
\tab return(S.homology(2))
\endcode
\medskip
\begincode
{\sl\# verify that the smallest example with nontrivial}
{\sl\# second homology is when $n = 143$}
for n in range(1,144):
\tab print(n, coprime\_free\_second\_homology(n))
\endcode

\goodbreak\bye